\documentclass[reqno,12pt]{article} \usepackage{fontenc} 
\usepackage{amsfonts}  \usepackage{amssymb}

\usepackage{pstricks} 

\def\alert#1{\smallskip{\hskip\parindent\vrule%
\vbox{\advance\hsize-2\parindent\hrule\smallskip\parindent.4\parindent%
\narrower\noindent#1\smallskip\hrule}\vrule\hfill}\smallskip}

\newcommand{\trdg}{\mbox{tr.deg.}}

\newcommand{\AAA}{\mathcal{A}}
\newcommand{\EEE}{\mathcal{E}}       
     
        \newcommand{\UUU}{\mathcal{U}}

    \newcommand{\htt}{\mbox{ht}}      \newcommand{\td}{\mbox{tr.deg.}}  
         
\newcommand{\RR}{\mathbb{R}}

\newcommand{\ifff}{\Longleftrightarrow}

\newcommand{\EE}{\mathbb{E}}        \newcommand{\ZZ}{\mathbb{Z}}

       \newcommand{\NN}{\mathbb{N}}
\newcommand{\xx}{\mathbf{x}} \newcommand{\ttt}{\mathbf{t}} 
\newcommand{\vv}{\mathbf{v}}
\newcommand{\rr}{\mathbf{r}}\newcommand{\bb}{\mathbf{b}}
\newcommand{\ee}{\mathbf{e}}\newcommand{\aaa}{\mathbf{a}}

\newtheorem{thm}{Theorem}[section]
\newtheorem{lem}[thm]{Lemma}
\newtheorem{cor}[thm]{Corollary}

\newtheorem{rem}[thm]{Remark}

\newtheorem{que}[thm]{Question}
\newtheorem{define}[thm]{Definition}
\usepackage{amssymb}
\usepackage{amsfonts}
\def\alert#1{\smallskip{\hskip\parindent\vrule%
\vbox{\advance\hsize-2\parindent\hrule\smallskip\parindent.4\parindent%
\narrower\noindent#1\smallskip\hrule}\vrule\hfill}\smallskip}

\usepackage{pstricks}

\begin{document}


\begin{center}
{\bf \Large Distances from the vertices of a regular simplex}\\~~\\
{Mowaffaq Hajja$^{1}$, Mostafa Hayajneh$^2$, Bach  Nguyen$^3$, and Shadi Shaqaqha$^4$}\\~\\
(1), (2), (4): Department of Mathematics, Yarmouk University, Irbid, Jordan\\
mowhajja@yahoo.com,~hayaj86@yahoo.com,~shadish2@yahoo.com\\~~\\
(3): Department of Mathematics, Louisiana State University, Baton Rouge, LA, USA\\
bnguy38@lsu.edu
\end{center}

\vspace{.4cm} \begin{quote}
{\small \bf Abstract.} {\small If $S$ is a given regular $d$-simplex of edge length $a$ in the $d$-dimensional Euclidean space $\EEE$, then the distances $t_1$, $\cdots$, $t_{d+1}$ of an arbitrary point in $\EEE$ to the vertices of $S$ are related by the elegant relation 
$$(d+1)\left( a^4+t_1^4+\cdots+t_{d+1}^4\right)=\left( a^2+t_1^2+\cdots+t_{d+1}^2\right)^2.$$ The purpose of this paper is to prove that this is essentially the only relation that exists
among $t_1,\cdots,t_{d+1}.$ The proof uses tools from analysis, algebra, and geometry.}
\end{quote}

\vspace{.4cm} \noindent {\bf Keywords:} algebraic dependence; Cayley-Menger determinant; germ;
height of an ideal; integral domain; Krull dimension;  Pompeiu's theorem; principal ideal;   real analytic function; 
regular simplex;      Soddy circles; transcendence degree

\vspace{.4cm} \noindent \section{Introduction} \label{int} Much has been written about the elegant relation
\begin{eqnarray}
(d+1) \left(a^4 + \sum_{j=1}^{d+1} t_j^4\right) = \left(a^2 + \sum_{j=1}^{d+1} t_j^2\right)^2 \label{rel}
\end{eqnarray}
that exists among the edge length $a$ of a regular $d$-dimensional simplex  in the Euclidean space $\RR^{d}$ and the distances $t_1,\cdots,t_{d+1}$ from the vertices 
of that  simplex  to an arbitrary point  in  $\RR^{d}$. The special case $d=2$ is illustrated in Figure 1 below, and the corresponding relation
\begin{eqnarray}
3 \left(a^4 + t_1^4 + t_2^4 + t_3^4 \right) = \left(a^2 + t_1^2 + t_2^2 + t_3^2 \right)^2 \label{rel-d=2}
\end{eqnarray}
was popularized by Martin Gardner in an article titled {\it Elegant triangle theorems not to be found in Euclid} that appeared in the June 1970 issue of {\it Scientific American}  \cite{MG-SA} and that was reproduced in \cite[Chapter 5, pp. 56--65]{MG}.
Feedbacks on Gardner's article appear in \cite{Gardiner},  \cite{Grimsey}, and \cite {Porteous}, and in possibly other places. A proof of the general case can be found in \cite{Bentin-1}, and another proof that uses  the Cayley-Menger formula for the volume of a simplex is given in \cite{prekites}. The relation (\ref{rel}) can also be derived using linear algebra.

 \vspace{.3cm}
One of the striking features of (\ref{rel}) is its symmetry, not only in $\{t_1,\cdots,t_{d+1} \}$ (which is only expected), but also in $\{ a, t_1,\cdots,t_{d+1} \}$. The fact that $a$ plays the same role in (\ref{rel}) as any other $t_i$ does not seem to have a satisfactory explanation\footnote{Recently, Dr. Ismail Hammoudeh of Amman University has come up with a satisfactory conceptual explanation of this symmetry. He will write this in a separate paper.}. Due to this symmetry, it is customary to set $a=t_0$ in (\ref{rel}) and write it in the form
\begin{eqnarray}
(d+1) \left(\sum_{j=0}^{d+1} t_j^4\right) = \left(\sum_{j=0}^{d+1} t_j^2\right)^2. \label{rel-2}
\end{eqnarray}

\vspace{0cm}
~~~~~~~~~~~~~~~~~~~~~~~~~\psset{unit=.3cm}
\begin{pspicture}(-11,-2)(11,20)
\psline[linestyle=dashed](-15,8)(-10,0)
\psline[linestyle=dashed](-15,8)(10,0)
\psline[linestyle=dashed](-15,8)(0,17.32)
\rput(-13.5,4){\small $t_2$} \rput(-2.5,3){\small $t_3$} \rput(-8.5,13.1){\small $t_1$}



\psline(10,0)(0,17.32)
\psline(-10,0)(0,17.32)
\psline(-10,0)(10,0)

\rput(0,-1){\small $a$} \rput(-6,8.7){\small $a$} \rput(6,8.7){\small $a$}
\rput(0,18){\small $\vv_1$} \rput(-10.5,-1){\small $\vv_2$} \rput(10.5,-1){\small $\vv_3$}
\rput(-15.5,8.2){\small $\ttt$}

\rput(0,-2.5){\small \bf Figure 1}
\end{pspicture}

\vspace{.5cm}
Another striking feature of (\ref{rel}) (or rather (\ref{rel-2}))  is its similarity with the relation
\begin{eqnarray}
d \left(\sum_{j=0}^{d+1} \left( \frac{1}{r_j}\right)^2\right) = \left(\sum_{j=0}^{d+1} \frac{1}{r_j}\right)^2
\label{rel2}
\end{eqnarray}
that exists among the {\it oriented} radii $r_0, \cdots, r_{d+1}$ of $d+2$ spheres in $\RR^{d}$ that are in mutual external touch.
The relation (\ref{rel2}) seems to have been known to Descartes, but it was rediscovered again and again by many. Its rediscovery by a Nobel prize winner, the chemist Frederick Soddy,  and the famous poem that he wrote to verbalize it, are probably major  reasons for its popularity. A proof of (\ref{rel2}) can be found, for example, in \cite{Pedoe} and \cite{Coxeter}, and
Soddy's poem {\it The Kiss Precise} first appeared in {\it Nature} in 1936 and  is reproduced in the last page of \cite{Oldknow} and in \cite[Chapter 3, pp.~33--34]{MG}. The relation (\ref{rel2}) is often referred to in the literature as {\it The Descartes' Circle Theorem}.

\bigskip It is now natural to ask whether there are other relations, beside the one in (\ref{rel}), among $t_1,\cdots,t_{d+1}$.
We state this  question precisely below, and we devote the rest of the paper to answering it.

\begin{que} \label{Q1}
{\em Let $S = [\vv_1,\cdots,\vv_{d+1}]$ be a regular $d$-simplex of side length $t_0$ in $\RR^d$. For every $\ttt \in \RR^d$,  let $t_j$ be the distance from $\ttt$ to $\vv_j$, $1\le j \le d+1$. Let 
$R = \RR[T_1,\cdots,T_{d+1}]$ 
be the ring of polynomials in the indeterminates   $T_1$, $\cdots$,  $T_{d+1}$ over the field $\RR$ of real numbers, and let 
$I$  be the ideal in $R$ defined by
\begin{eqnarray}
I&=& \{ f \in R : f\left(t_1,\cdots, t_{d+1}\right) = 0 ~~\forall~ \ttt \in \RR^{d}\}, \label{I}
\end{eqnarray}
and let  \begin{eqnarray}
F &=& (d+1) \left( t_0^4 + \sum_{j=1}^{d+1} T_j^4\right) - \left( t_0^2 + \sum_{j=1}^{d+1} T_j^2\right)^2. \label{F}
\end{eqnarray}
Is $I$  generated by $F$?}  
\end{que}

\bigskip
We shall establish  an affirmative answer to this  question in Theorem \ref{main} in Section \ref{MAIN}.  Readers may find it interesting that the ideas used in the proof of this innocent-looking theorem are so varied, involving tools from analysis, algebra, and geometry. As some of this material cannot be assumed to be known to the potential readers of this article, we have chosen to be self-contained, writing the definitions of the terms, and giving adequate references to, or proofs of,  the theorems used.
 Readers are advised to read the main theorem, Theorem \ref{main}, and its proof first, and decide for themselves what sections of the paper they need to go back to and read.

The paper is organized as follows. 
Section \ref{CM} recalls a theorem, proved in \cite{Van}, pertaining to the irreducibility of a class of Cayley-Menger determinants. 
Section \ref{RAF}  introduces the preliminary definitions and theorems from the theory of real analytic functions. The only  reference that we have used is the book \cite{RAF}.  Section \ref{AI} puts together the necessary algebraic tools.  These include the height of an ideal, algebraic independence over a field, transcendence bases and degree of an extension, the Krull dimension of a ring, and how these are related. The only reference that we referred to here is R. Y. Sharp's book  \cite{Sharp}. We have also included in this section 
that the distances from an arbitrary point in $\RR^n$ to any $n$ vertices of a regular $n$-simplex in $\RR^n$ are algebraically independent (over $\RR$).   Interestingly, an essential step in the proof follows from a simple fact that happened to appear in Euclid's {\it Elements}, namely as Propositions 7 and 8 of Book III.  Section \ref{MAIN} contains the main theorem and its proof.  The last section, Section 
\ref{Q}, contains  a list of problems that may generate  further research. We expect that most, or all, of these problems are within reach of a young researcher, and we also expect the material in Sections  \ref{CM}, \ref{RAF}, \ref{AI}, and \ref{MAIN} to be useful to such a researcher and to others, and easier to refer to than to refer to various books.

\section{Irreducibility of a class of Cayley-Menger determinants}\label{CM}
The following theorem is taken from \cite{Van}. Its special case $t = n = d+1$, $d \ge 2$ shows that the polynomial 
$$F = (d+1)\left(a^4 + \sum_{j=1}^{d+1} T_j^4\right) -
\left(a^2 + \sum_{j=1}^{d+1} T_j^2\right)^2$$ that appears in Question \ref{Q1} is irreducible. This fact will, in turn, be  used in the proof of the main theorem.

\begin{thm} \label{Van}
Let $F = F(x_1,\cdots,x_n) \in \RR[x_1,\cdots,x_n]$ be defined by
\begin{eqnarray}\label{FFF}
F &=& t\left(a^4 + \sum_{j=1}^{n} x_j^4\right) -
\left(a^2 + \sum_{j=1}^{n} x_j^2\right)^2,
\end{eqnarray}
and suppose that $t \ne 0$. Then $F$ is irreducible except when $n \le 2$ and when $(n,t) = (3,2)$.
\end{thm}

\vspace{.1cm} \noindent {\it Proof.} This is taken from \cite{Van}. \hfill $\Box$

\vspace{.5cm}
The situations when the field $\RR$, appearing in the previous theorem, is replaced by an arbitrary field, and when $n$ is not restricted to the values $n \ge 3$ are treated in full detail in the aforementioned reference \cite{Van}. The relation of the polynomial in (\ref{FFF}) to the Cayley-Menger determinant is exhibited in  \cite{prekites}.


\section{Real analytic functions in  several variables}\label{RAF}
In this section, we present the basic material on real analytic functions that will be needed in the proof of the main theorem, Theorem \ref{main}.  
The treatment is self contained, and all necessary definitions are given. We feel that this  interesting subject (of real analytic functions) is not usually covered in standard required courses in graduate schools, and we also feel that there are not many textbooks on the subject. Our only reference is \cite{RAF}, and any differences between our presentation and that in \cite{RAF} are very slight and trivial, and they are made with the permission of the first author of \cite{RAF}.

\begin{define}\label{DEF:ps}
{\em  (\cite[Definition 2.1.4,  p. 27]{RAF}) 
Let $\ZZ^+ = \NN \cup \{0\}$ be the set of non-negative integers, and let $\left(\ZZ^+ \right)^m$, $m \in \NN$, be denoted by $\Lambda (m)$. If $\ee = (e_1,\cdots,e_m) \in \Lambda (m)$, and if $\rr = (r_1,\cdots,r_m) \in G^m$, where $G$ is any commutative ring with 1, then $\rr^{\ee}$ stands for the product  $$\prod_{j=1}^m r_j^{e_j}.$$
A power series in the $m$ variables $\xx = (x_1,\cdots,x_m)$ with center at $\aaa = (a_1,\cdots, a_m) \in \RR^m$ is a formal expression  
\begin{eqnarray}
\sum_{\ee \in \Lambda (m)} c_{\ee}  (\xx-\aaa)^{\ee}, \mbox{~where $c_{\ee} \in \RR.$} \label{PS}
\end{eqnarray}
The power series (\ref{PS}) is said to converge at $\bb = (b_1,\cdots,b_m) \in \RR^m$ if some rearrangement of it converges. More precisely, the power series (\ref{PS}) converges at $\bb = (b_1,\cdots,b_m)$ if there is a bijection  $\phi : \ZZ^+ \to \Lambda (m)$ such that the sequence of partial sums of the series $$\sum_{j=0}^{\infty} c_{\phi (j)} (\bb - \aaa)^{\phi (j)}$$ converges.}
\end{define}

\begin{rem}\label{REM:contrast}
{\em In first courses of calculus, a series $\sum t_j =\sum_{j=1}^{\infty} t_j$ of real (or complex) numbers was said to converge if and only if its sequence of partial sums converges. Thus if $c_j = (-1)^j / j$ for $j \in \NN$, then the sequence of partial sums of the series $\sum c_j$ converges. Now let  $\sum d_j$ be a rearrangement of $\sum c_j$ such that the sequence of partial sums of the series $\sum d_j$ does not converge, and consider the power series $\sum d_j x^j$. According to the definition above, this series converges at $x=1$; but according to the definition given in calculus books, the series does not converge at $x=1$. Thus the definition in caculus books and the definition above give two different {\it sets} $S_{old}$ and $S_{new}$ of convergence of $\sum d_nx^n$, with $1 \in S_{new}$ and $1 \notin S_{old}$.  However, we shall see that  these two sets have the same interior. This interior, usually referred to as the {\it domain} of convergence,  is what really matters for us.

Note also that if a series $\sum t_j$ converges, then it does not necessarily have a well-defined sum. This is because the series may have two different convergent rearrangements having different sums. In fact, if a series of real numbers is convergent but not absolutely convergent, then it has a rearrangement that converges to any prescribed number; see Theorem 12-33 (p.~368) of \cite{Apostol-MA}. However, if  the series $\sum t_n$ converges absolutely, in the sense that $\sum |t_n|$ converges,
then both senses of convegence coincide, and the  series will have a well-defined sum; see Theorem 12-32 (p.~367) of \cite{Apostol-MA}.
In particular, this holds for convergent series whose terms are non-negative.


We shall also see below that if a power series converges at every point in an open set $U$, then it converges absolutely  at every point in $U$, and thus it defines a function on $U$.}
\end{rem}

\begin{thm} \label{THM:Abel} {(\bf  \cite[Proposition 2.1.7 (Abel's Lemma),  p. 27]{RAF})} If the power series $\sum_{\ee \in \Lambda (m)} c_{\ee} \xx^{\ee}$ converges at a point $\xx=\bb = (b_1, \cdots, b_m) \in \RR^m$, then it converges uniformly and absolutely on compact subsets of the silhouette of $\bb$, i.e., the open box $$\left( -|b_1|,|b_1| \right) \times \cdots \times \left( -|b_m|,|b_m| \right).$$
\end{thm}

\begin{cor} \label{COR:c->ac}  If a power series converges at every point in an open set $U$, then it converges absolutely at every point in $U$.
Consequently, it defines a function on $U$.\end{cor}

\vspace{.09cm} \noindent  {\it  Proof.} Let $q \in U$. Then $U$ contains a closed  box $B$ centered at $q$. Thus $q$ belongs to the silhouette $S$ of  some vertex, say $\vv$, of $B$. Since  the given power series converges at $\vv$, it follows  from Abel's lemma that it converges absolutely on $S$, and hence at $q$. \hfill $\Box$

\begin{define} \label{DEF:ra}   {\em {(\bf  \cite[Definition 2.2.1, p. 29]{RAF})} We say that  the real-valued function $f$ is {\it real analytic} at a point $p$  in  $\RR^n$ if   $p$ has a  neighborhood $U$ on which  $f$ can be represented  as a power series centered at $p$. We say that  $f$ is real analytic on a set $U$ if $f$ is real analytic at every point of $U$.}
\end{define}


\begin{thm}  \label{THM:open}
If $f$ is real analytic at $p \in \RR^n$, then  $f$ is real analytic on a neighborhood $U$ of $p$. Consequently, the set where $f$ is real analytic is open. \end{thm}

\vspace{.1cm} \noindent {\it Proof.} Let $F_p(z)$ be a power series centered at $p$ (i.e., in powers of  $(z-p)$) that represents $f$ on some neighborhood $U$ of $p$. By Proposition 2.2.7 (p. 32) of \cite{RAF}, $F_p(z)$ is real analytic on $U$. But $F_p(z)$ coincides with $f$ on $U$. Therefore $f$ is real analytic on $U$. \hfill $\Box$

\begin{thm} \label{THM:ring} {\bf  (\cite[Theorem 2.2.2, p. 29]{RAF})}
If $f$, $g$ are real analytic on the  subsets $U$, $V$ of $\RR^n$, respectively, then $f+g$, $fg$ are real analytic on $U \cap V$, and $f/g$ is real analytic on $U \cap V \cap \{p : g(p) \ne 0 \}$.
\end{thm}

\begin{rem} \label{REM:pol=entire}{\em Needless to say, polynomials are {\it real entire} (= analytic at all points).} 
\end{rem}

\begin{thm}  \label{THM:f=g} {\bf  (\cite[Corollary  1.2.7, p. 14]{RAF})} If   $f$ and $g$ are real analytic on an open interval $U \subseteq \RR$ and if there is a sequence of distinct points $x_1$, $x_2$, $\cdots$ in $U$ with $x_0 =\lim_{n\to \infty} x_n \in U$ and such that $f(x_n) = g(x_n)$ for all $n \in \NN$, then $f(x) = g(x)$ for all $x \in U$. \end{thm}

\begin{cor}  \label{COR:acc-pts} 
{\bf  (\cite[Corollary  1.2.7, p. 14]{RAF})} 
If   $f$ is real analytic on the open interval $U \subseteq \RR$, and if  the set of zeros of $f$ has an accumulation point in $U$, then $f$ is identically zero on $U$. \end{cor}

\begin{thm} \label{THM:ID-1}
The set of real analytic functions on an open interval $U \subseteq \RR$ is an integral domain. \end{thm}

\vspace{.1cm} \noindent {\it Proof.} Let $f$ and $g$ be real analytic on an open interval $U$, and suppose that $fg$ is identically zero on $U$. We are to show that either $f$ or $g$ is identically zero on $U$.

Let $K$ be a compact subinterval of $U$, and let $Z(f)$ and $Z(g)$ be the sets of zeros of $f$ and $g$, respectively, in $K$. Then $Z(f) \cup Z(g) = K$, and therefore either $Z(f)$ or $Z(g)$ is an infinite set. Suppose that $Z(f)$ is an infinite set. By the Heine-Borel theorem,  $Z(f)$ has an accumalation point in $K$. Therefore the set of zeros of $f$ in $U$ has an accumalation point in $U$. By Corollary \ref{COR:acc-pts}, 
$f$ is identically zero on $U$. This proves our claim.
\hfill $\Box$

\begin{rem} \label{f=xy} {\em Theorems \ref{THM:f=g} and \ref{COR:acc-pts} do not seem to have  analogues in higher dimensions. For example, the function $f : \RR^2 \to \RR$ defined by $f(x,y) = xy$ is real entire, and its zero set consists of the $x$- and $y$- axes. But $f$ is not identically zero. However, Theorem \ref{THM:ID-1} does have an analogue in higher dimensions; see Theorem \ref{THM:ID-2} and Corollary \ref{xxx}.}
 \end{rem}

\begin{lem} \label{LEM:x->a}
Suppose that $f = f(x_1,\cdots,x_n)$ is real analytic at $a = (a_1,\cdots,a_n)$, and let $1 \le k \le n$. Let $$f^* = f^* (x_{k+1},\cdots,x_n) = f(a_1,\cdots, a_k, x_{k+1},\cdots,x_n).$$ Then $f^*$ is real analytic at $a^* = (a_{k+1},\cdots,a_n)$. The same holds for the function $f_*$ and the point $a_*$ defined by
$$f_* = f_* (x_{1},\cdots,x_k) = f(x_1,\cdots, x_k, a_{k+1},\cdots,a_n),~~a_* = (a_1,\cdots,a_k).$$
\end{lem}

\vspace{.1cm} \noindent {\it Proof.} It is clearly enough to prove the case $k=1$, since the general case follow by repeated application.

Let $F = F(x_1,\cdots,x_n)$ be a power series centered at $a$ and represents $f$ on a neighborhood $U$ of $a$. Let $V$ be the hyperplane defined by $x_1 = a_1$, and let $U^*$ be the projection of $U \cap V$
on the $(x_2,\cdots,x_n)-$hyperplane. Thus
$$(a_2,\cdots,a_n) \in U^* \ifff (a_1,a_2,\cdots,a_n) \in U \cap V.$$
Then $U^*$, being the image of the open set $U \cap V$ under the projection $(x_1,x_2,\cdots,x_n) \mapsto  (x_2,\cdots,x_n)$, is a neighborhood of $a^*$. Also, the power series $F^* = F^* (x_2,\cdots,x_n) $ defined by $$F^* (x_2,\cdots,x_n) = F(a_1,x_2,\cdots,x_n)$$
represents $f^*$ on $U^*$. Therefore $f^*$ is real analytic at $a^*$, as desired. \hfill $\Box$

\begin{thm}  \label{THM:ID-2}  Let $U_i$, $1 \le i \le n$, be  open  intervals in $\RR$, and let $U$ be the open box  in $\RR^n$ defined by $U = U_1 \times \cdots \times U_n$. Then the set of functions that are real analytic on $U$ is an integral domain. \end{thm}

\vspace{.1cm} \noindent {\it Proof.} Let $f, g$ be real analytic on $U$, and suppose that $fg$ is identically zero on $U$. We claim that  either $f$ or $g$ is identically zero on $U$.

By Theorem \ref{THM:ID-1}, the claim is   true when $n=1$. So we proceed by induction.

Let $a \in U_n$, and let $f_a$ and $g_a$ be defined on $V = U_1 \times \cdots \times U_{n-1}$ by
$$f_a (x_1,\cdots,x_{n-1}) = f (x_1,\cdots,x_{n-1},a),~
g_a (x_1,\cdots,x_{n-1}) = g (x_1,\cdots,x_{n-1},a).$$
Then $f_a$ and $g_a$ are real analytic on $V$, by Lemma \ref{LEM:x->a}, and $f_ag_a$ is zero there. This is because
\begin{eqnarray*}
f_ag_a (x_1,\cdots,x_{n-1}) &=&  f (x_1,\cdots,x_{n-1},a) g (x_1,\cdots,x_{n-1},a)\\
&=& fg (x_1,\cdots,x_{n-1},a) ~=~ 0.
\end{eqnarray*}
By the inductive assumption,  either $f_a$ or $g_a$ is zero on $V$. This is true for all $a$ in $U_n$. Therefore either $f_a$ or $g_a$ is zero on $V$ for infinitely many values of $a$ in $U_n$.
Without loss of generality, we may assume that 
\begin{eqnarray} \label{added}
\mbox{$f_a$ is zero on $V$ for infinitely many values of $a$ in $U_n$}.
\end{eqnarray}
 Now take any $p=(a_1,\cdots,a_{n-1})$ in $V$. The function $h$ defined on $U_n$ by $$h (x_n) = f (a_1,\cdots,a_{n-1},x_n)$$
is real analytic on $U_n$, by Lemma \ref{LEM:x->a}. Also, $h$ has infinitely many zeros in  $U_n$, by (\ref{added}). Therefore $h$ is zero on $U_n$, by Corollary \ref{COR:acc-pts}. Thus
$$f (a_1,\cdots,a_{n-1},x_n) = 0$$ for all $x_n \in U_n$. Since $(a_1,\cdots,a_{n-1})$ was an arbitrary point in $V$, it follows that
$$f (x_1,\cdots,x_{n-1},x_n) = 0$$
for all $(x_1,\cdots,x_{n-1},x_n) \in U$. Thus $f$ is zero on $U$, as desired. \hfill $\Box$

\begin{cor} \label{xxx}  Let $V$  be  an open set  in $\RR^n$. Then the set of functions that are real analytic on $V$ is an integral domain. \end{cor} 

\vspace{.1cm} \noindent {\it Proof.} Let $U$ be an open box contained in $V$. Let $S_V$ and $S_U$ be the rings of functions that are real analytic on $V$ and $U$, respectively. Also, $S_V$ is a subring of $S_U$, and $S_U$ is an integral domain. Therefore $S_V$ is an integral domain. \hfill $\Box$

\begin{rem} 
{\em On the set of all functions that are real analytic at a point $p \in \RR^m$, we define an equivalence relation $\equiv$ by $$f \equiv g \ifff f = g \mbox{~on some open neighborhood $U=U_{f,g}$ of $p$}.$$
Each equivalence class is called a {\it germ} or a $p$-germ. It is not difficult to define addition and multiplication  on the set $G$ of $p$-germs so that  $G$ becomes a ring. Now the proof above can  be mimicked to show that $G$ is an integral domain. This is stronger than Corollary \ref{xxx}. } \end{rem}

\begin{thm} \label{aaaa}{\bf ( \cite[Proposition 2.2.8, p. 33]{RAF})}
If $f_j$, $1 \le j \le d$, are real analytic at $p \in \RR^n$, and if $g$ is real analytic at $(f_1(p),\cdots,f_d(p)) \in  \RR^d$, then the composition $g(f_1,\cdots,f_d)$ is real analytic at $p \in  \RR^n$. \end{thm}

\begin{thm} \label{bbbb} ({\bf \cite[Theorem 2.3.1,  p. 35, and Theorem 2.3.5, p. 40]{RAF})} Suppose that $f(x_1,\cdots,x_n;y)$ is real analytic at a point $(p_1,\cdots,p_n;q) \in \RR^n \times \RR = \RR^{n+1}$, and suppose that 
$f(p_1,\cdots,p_n;q) = 0$ and that 
$$\frac{\partial f}{\partial y} (p_1,\cdots,p_n;q) \ne 0.$$ Then there exists a  neighborhood $V$ of $(p_1,\cdots,p_n)$ and a function $\Phi : V \to \RR$ that is real analytic at $(p_1,\cdots,p_n) \in \RR^n$ such that
$\Phi (p_1,\cdots,p_n) = q$ and such that $$f(t_1,\cdots,t_n;\Phi (t_1,\cdots,t_n)) = 0$$
for all $(t_1,\cdots,t_n)$ in some  neighborhood of $(p_1,\cdots,p_n)$.
\end{thm}

\begin{cor} \label{yyy}
If $f$ is  real analytic on some open set $U \subseteq \RR^n$, and if $f(p) > 0$ for some $p \in U$, then there exists a neighborhood $W \subseteq U$ of $p$ such that $\sqrt{f}$ exists and is  real analytic on $W$.  In particular, if $p \in \RR^n$, then the function $h : \RR^n \to \RR$  defined by $h(x) = \| x - p\|$ is  real analytic at all points except  at $p$.
\end{cor}

\vspace{.1cm} \noindent {\it Proof.} Suppose that $f$ is real analytic on some open set $U \subseteq \RR^n$, and that $f(p) > 0$ for some $p \in U$. Let  $V =  U \times \RR$, and let $g : V \to \RR$ be defined by $$g(x,y) = f(x) - y^2, ~~x \in U,~y\in \RR.$$
Since $f$ is real analytic on $U$, it is real analytic on $U \times \RR$. Similarly, $y^2$ is real analytic on $\RR$. So it is real analytic on $\RR^n \times \RR$ and hence on $V$. Therefore $g$ is real analytic on $V$. Let $q = \sqrt{f(p)}$. Then $g$ is real analytic at $(p,q)$. Also, $g(p,q) =0$, and $$\frac{\partial{g}}{\partial{y}}\left|_{(p,q)} \right.= - 2y\left|_{(p,q)} \right.= - 2q \ne 0.$$ Therefore there exist, by Theorem \ref{bbbb}, a neighborhood $W \subseteq U$ of $p$ and a real analytic function  $\Phi$ on $W$ such that $f(x) = [\Phi (x)]^2$ on $W$, and $W$ can be chosen such that $f(x) > 0$ for all $x$ in $W$. Since $f > 0$ on $W$, it follows that $\Phi$ is never 0 on $W$, and therefore either $\Phi$ is positive on $W$ or 
$\Phi$ is negative on $W$. In the latter case, we replace $\Phi$ by $- \Phi$, and obtain
$$\Phi > 0 \mbox{~on~} W,~~f(x) = [\Phi (x)]^2  \mbox{~on~} W.$$
Thus $\Phi = \sqrt{f}.$

To prove the last statement, 
let $g(x) =[h(x)]^2 = (x_1-p_1)^2 + \cdots + (x_n-p_n)^2.$
It is clear that the partial derivatives of $h = \sqrt{g}$ at $x=p$ do not exist. Therefore $h = \sqrt{g}$ is not real analytic at $p$. Now let $q=(q_1,\cdots,q_n)\in \RR^2$ be any point other than $p$. Then $g$, being a polynomial, is real analytic everywhere, and $g(q) > 0$. Therefore
$\sqrt{g}$ exists on a neighborhood of $q$ and is  real analytic there. It is clear that  $\sqrt{g}$ coincides with $f$.
\hfill $\Box$

\section{Algebraic independence of certain functions} \label{AI}
Let $\{\vv_1,\cdots,\vv_d\}$ be a set of affinely independent points in the $d$-dimensional euclidean space $\RR^d$, and let  the functions $\phi_j : \RR^d \to \RR$, $1 \le j \le d$, be defined  by $\phi_j(\xx) = \| \xx - \vv_j\|.$ Let $\UUU$ be an everywhere dense subset of $\RR^d$, and let $\AAA$ be the ring of real-valued functions that are real analytic on $\UUU$. We have seen in Corollary \ref{xxx}  that $\AAA$ is an integral domain, and we have seen in Corollary \ref{yyy} that the functions $\phi_j$, $1 \le j \le d$, belong to $\AAA$. In this section, we shall  show that these functions are algebraically independent over $\RR$.

We start by the necessary  preliminaries that we need from algebra.

\subsection{Algebraic preliminaries.}
Our reference is \cite{Sharp}, and {\it all rings referred to are commutative with 1}. If $B$ is any ring, then
$B[T_1,\cdots,T_n]$ stands for  the polynomial ring over $B$ in the $n$ indeterminates $T_1,\cdots,T_n$.

Let $S$ be a ring (commutative with identity $1 \ne 0$).

If $P$ is a prime ideal in $S$, then the {\it height} of $P$, written $\htt P$, is the supremum of all non-negative integers $k$ for which there exists a sequence $P_0, P_1,\cdots,  P_k = P$ of prime ideals in $S$ of the form 
$$P_0 \subset P_1 \subset \cdots \subset P_k = P,$$ where $\subset$ stands for {\it strict} inclusion. 
If the supremum does not exist, we write $\htt P = \infty$. The {\it dimension} of $S$,   written $\dim S$, is the supremum of all non-negative integers $k$ for which there exists a sequence $P_0, P_1,\cdots,  P_k = P$ of prime ideals in $S$ of the form 
$$P_0 \subset P_1 \subset \cdots \subset P_k = P.$$
If the supremum does not exist, we write $\dim S = \infty$. 
Thus $$\dim S = \sup \{ \htt P : P \mbox{~is a prime ideal in $S$}\},$$
where the supremum is defined for all subsets of $[- \infty, +\infty]$. It is easy to see, as done in Remark 14.18 (viii) (page 279) of \cite{Sharp},  that
\begin{eqnarray}
\htt~ P + \dim S/P &\le& \dim S, \label{main-1}
\end{eqnarray}
where we adopt the convention  that $\infty + \infty = \infty$, and 
$\infty + n = \infty$ for all integers $n$.

Now let $R$ be a subring of $S$. 



A finite subset $\{ s_1,\cdots,s_n\}$ of $S$ is said to be {\it algebraically independent} over $R$ if there does not exist a non-zero polynomial $f = f(T_1,\cdots,T_n)$ in the polynomial ring $R[T_1,\cdots,T_n]$ such that $f(s_1,\cdots,s_n) = 0.$  An arbitrary subset of $S$ is said to be {\it algebraically independent} over $R$ if every  finite subset of $S$  is algebraically  independent over $R$. This is Definition 1.14 (page 8) of \cite{Sharp}. If both $R$ and $S$ are fields, then a subset $B=\{ s_1,\cdots,s_n\}$ of $S$ that is algebraically independent over $R$ is called a (finite) {\it  transcencdence basis} of $S$ over $R$ if every subset of $S$ that properly contains $B$ is algebraically dependent over $R$. This is Definition 12.54 (page 239) of \cite{Sharp}. By Theorem 12.53 (page 239) of \cite{Sharp}, if $S$ has a finite transcendence basis over $R$, then any two such  bases have the same number of elements. This number is called the {\it  transcendence degree} of $S$ over $R$, and is denoted by $\trdg_R S$.

If $S$ is an integral domain that contains a field $R$ and that is finitely generated over $R$, i.e., $S$ is an affine $R$-algebra, and if $L$ is the field of quotients of $S$, then
\begin{eqnarray}
\dim S &=& \trdg_R L. \label{main-2}
\end{eqnarray}
This is Corollary 14.29 (page 282) in \cite{Sharp}.
It follows that if $R$ is a field, 
then 
\begin{eqnarray}
\dim  R[T_1,\cdots,T_n]  &=& n. \label{main-3}
\end{eqnarray}

We also shall need the facts that if $R$ is a field (or any unique factorization domain), then the polynomial ring $R[T_1,\cdots,T_n]$ is a unique factorization domain, and that an irreducible element in a unique factorization domain $D$ generates a prime ideal of $D$;   see \cite[Theorem 1.42, p.~17]{Sharp} and \cite[Exercise 3.42, p.~49]{Sharp}. Now the following theorem, that we will use later, follows immediately.

\begin{thm} \label{ht=1} Let $R$ be a field,  and let 
$P$ be a prime ideal of $R[T_1,\cdots,T_n]$ that contains an irreducible polynomail $f$. If $\htt~ P = 1$, then $P$ is generated by $f$. \end{thm}

\vspace{.2cm} \noindent {\it Proof.} Let $Q$ be the ideal generated by $f$. Then $Q$ is a prime ideal. If $Q \ne P$, then the chain $P \supset Q \supset \{0\}$ would imply the contradiction that  $\htt P \ge 2$. Therefore $P=Q$, as claimed. \hfill $\Box$

\subsection{Algebraic independence of distances to affinely independent points.}
In Theorem \ref{zzz} below, we prove that the functions $\phi_j : \UUU \to \RR$, $1 \le j \le d$, defined in the previous section  are algebraically independent over $\RR$. To do so, we need the  simple lemma, Lemma \ref{III78}. 
Interestingly, this lemma is essentially nothing but a combined form of Propositions 7 and 8 of Book III of Euclid's {\it Elements}. These  propositions, illustrated in Figure 2 below,  state that if $UV$ is a diameter of a circle centered at $A$, and if $B$ lies in the open line segment $UA$, then as a point $P$ moves from $U$ to $V$ along one of the semicircles having diameter $UV$, the length of the line segment $BP$ increases strictly. This is clearly an immediate consequence of Proposition 24 of Book I of Euclid's {\it Elements}, better known as {\it the open mouth theorem.}

\begin{lem} \label{III78} Let $\EE$ be a Euclidean space of dimension greater than or equal to 2. Let $\Gamma$ be a circle in $\EE$ and let $A$ be its center. Let  $D \in \EE$, and suppose that the orthogonal projection $B$ of $D$ on the plane $H$ of $\Gamma$ is not $A$.
Then  as $P$ ranges in $\Gamma$, $\| D-P\|$ assumes infinitely many values.
\end{lem}

\vspace{.2cm} \noindent {\it Proof.} 
For $P \in \Gamma$, we have
$$\|D - P\|^2 = \|D - B\|^2 + \|B - P\|^2.$$
Since $\| D-B\|$ is fixed, it is sufficient to prove our statement for the function $g(P) := \|B - P\|$.
We also may assume that $H$ is the usual Euclidean plane.

Let $U$ and $V$ be the points where the line $AB$ meets $\Gamma$. As  $P$ moves from $U$ to $V$ on one of the semicircles of $\Gamma$, it follows from Propositions 7 and 8 of Book III of Euclid's {\it Elements} that the quantity $\| B - P\|$ changes strictly monotonically from $\| B - U\|$ to $\| B - V\|$. Since $\| B - U\| \ne \| B - V\|$, it follows that $\|B - P\|$ assumes infinitely many values. \hfill $\Box$

\begin{center}
\psset{unit=.45cm}
\begin{pspicture}(-6,-6)(7,6)
\psarc(0,0){5}{0}{180}  \psarc[linestyle=dashed](0,0){5}{180}{360}
\rput(-5.6,-.0){\small $U$} \put(5.25,-.25) {\small $V$}
 \put(-4.00,-.9){\small $B$} \rput(0.00,-.7){\small $A$} 
\put(-1.5,5.1){\small $P$} \put(4.2,3.1){\small $P$}  \put(-3.7,4.1){\small $P$}

\psline(-3.6,0)(-3,4) \psline(-3.6,0)(4,3) 
\psline(-3.6,0)(-1.294095224,4.829629132)
\psline(-5,0)(5,0)

\rput(0,0){$\bullet$}  \rput(-3.6,0){$\bullet$}
\rput(0,-6.7){\bf Figure 2}


\end{pspicture}
\end{center}

\vspace{.5cm} \noindent  {\bf Note.}
Proposition 7 of Book III of Euclid's {\it Elements} deals with the case when $B$ lies inside $\Gamma$, while Proposition 8 deals with the case when $B$ lies outside $\Gamma$.

\begin{thm} \label{zzz}
Let $n \ge 1$, and let $\EE$ be a Euclidean space of dimension  greater than or equal to  $n$, and let 
$\vv_1,\cdots,\vv_{n}$ be affinely independent points in $\EE$. 
Let the functions $\phi_1,\cdots,\phi_n : \EE \to \RR$ be defined by
$$\phi_j  (p) = \|p - \vv_j\|,~ 1 \le j \le n.$$
Then the functions $\phi_1,\cdots,\phi_n$, (thought of as elements in the ring  of all real-valued functions on $\EE$) are algebraically independent (over $\RR$). 
\end{thm}

\vspace{.1cm} \noindent {\it Proof.}
The  claim is trivial for $n =1$. In fact, if $\vv_1=\vv$ is any (affinely independent) point in a Euclidean space of dimension at least 1,  if $\phi :  \EE \to \RR$ is defined by $\phi (p) = \| \vv - p\|$,  and if $f(T)$ is a non-zero polynomial in $\RR[T]$ such that $f(\phi) = 0$, then $f(\phi (p)) = 0$ for all $p$ in $\EE$. Since $\{ \| p - \vv \|: p \in \EE\}$ is an infinite set, we obtain the contradiction that $f$ has infinitely many zeros.

We now prove our claim for $n=2$, since this is the step that we will use in our proof by induction.  Thus let  $\EE$ be Euclidean space of dimension at least 2, and let $\vv_1=A$ and $\vv_2=B$ be points in $\EE$ that are affinely independent,  i.e., $A\ne B$. Let $\alpha$ and $\beta$ be the functions from $\EE$ to $\RR$  defined by $\alpha (p) =  \| p - A\|$ and $\beta (p)  =   \| p - B\|$. 
We are to show that $\alpha$ and $\beta$ are algebraically independent over $\RR$.
Since $\dim_{\RR} \EE \ge 2$, it follows that $\EE$ has a 2-dimensional subspace $H$ that contains $A$ and $B$. We  work within $H$, and thus we may assume that $\EE = H$, i.e., the ordinary Euclidean plane.

Suppose that there is a non-zero polynomial $f = f(X,Y)$ in the polynomial ring $\RR[X,Y]$ such that $f(\alpha, \beta) = 0$, i.e., $f(\alpha (p), \beta  (p)) = 0$ for all $p$ in $\EE$. If $f$ is a constant polynomial, then it would have to be the zero polynomial. So we may assume that one of the variables $X$ and $Y$, say $Y$ appears in $f$.
Thus $f$ can be written in  the form $$f(X,Y) = Y^d g_d(X) + Y^{d-1}  g_{d-1} (X) + \cdots + Y g_1(X) + g_0 (X),$$ where $d \ge 1$, $g_j \in \RR[X]$, and $g_d$ is not the zero polynomial. Thus there exists a positive number $r$ such that $g_d (r) \ne  0$. 
Let $\Gamma$ be the circle centered at $A$ and having radius $r$. Then for all $p \in \Gamma$, $$g_d (\alpha (p)) = g_d (\| A - p\|) = g_d (r) \ne 0.$$ Since $A\ne B$, the line joining $A$ and $B$ crosses $\Gamma$ at the end points $U$, $V$ of a diameter. By Lemma \ref{III78}, the set $\{\beta (p) : p \in \Gamma\}$ is infinite. For $p \in \Gamma$,  $f(\alpha (p), \beta (p)) = f(r, \beta (p)) = 0$, since $f(\alpha, \beta) = 0$. Therefore the polynomial $h(Y):= f(r, Y)$ has infinitely many zero. This contradicts the fact that
$h(Y)$ is non-zero since the coefficient of $Y^d$ in $h(Y)$ is  $g_d(r) \ne 0$. This proves the case $n=2$.

We now proceed by induction. Thus we suppose that the claim is true for $n=k$ for some $k \ge 2$, and we  are to show that the claim is true for $n=k+1$. So we let $\vv_1,\cdots,\vv_{k+1}$ be given affinely independent points in some Euclidean space $\EE$ having dimension at least $k+1$, and we define the functions $\phi_j : \EE \to \RR$, $1 \le j \le k+1$, as before, i.e., by
$$\phi_j  (p) = \|p - \vv_j\|,~ 1 \le j \le k+1.$$
We are to show that $\phi_j$, $1 \le j \le k+1$, are algebraically independent (over $\RR$).

We may clearly assume that $\dim_{\RR} \EE = k+1$.
This is because if certain  functions defined on a set are algebraically dependent, then their restrictions to any subset remain algebraically dependent.

Suppose that $\phi_j$, $1 \le j \le k+1$, are algebraically dependent (over $\RR$).
Then there exists a non-zero polynomial $g$ in $\RR[T_1,\cdots,T_{k+1}]$ such that $g(\phi_1,\cdots,\phi_{k+1}) = 0$, i.e., 
\begin{eqnarray} \label{g-phi}
g(\phi_1 (p), \cdots, \phi_{k+1}  (p)) = 0 ~~\forall ~p \in \EE.
\end{eqnarray}
Clearly $g$ cannot be the constant 
polynomial, and hence one of the variables $T_1, \cdots, T_{k+1}$ must appear in $g$. Without loss of generality, we may assume that $T_{k+1}$ appears in $g$.
Therefore  $g$ can be written in the form
$$g = g_d T_{k+1}^d + \cdots + g_1 T_{k+1} + g_0,$$
where $d \ge 1$, where $g_j$, $1 \le j \le d$, belong to $\RR[T_1,\cdots,T_k]$, and where $g_d$ is not the zero polynomial. Since $\phi_1, \cdots, \phi_k$ are algebraically independent, it follows that $g_d(\phi_1,\cdots,\phi_k)$ is not the zero function on $\EE$. Therefore there exists a point $p \in \EE$ such that $g_d(\phi_1(p),\cdots,\phi_k(p)) \ne 0$. By the continuity of $\phi_i$'s and $g_d$, there exists an open set $U$ in $\EE$ such that \begin{eqnarray}\label{STAR}
g_d(\phi_1(q),\cdots,\phi_k(q)) \mbox{~ is nonzero for every  $q \in U$}.
\end{eqnarray}
We choose $q \in U$ such that $q$ does not lie in the affine hull $H_{k+1}$ of $\vv_1, \cdots, \vv_{k+1}$. This is possible since $U$ is open in $\EE$ and since 
$H_{k+1}$, as a subset of $\EE$, has an empty interior. Again this is so because 
$\dim H_{k+1} = k   < k+1 = \dim \EE$.

Let $H_k$ be the affine hull of
$\vv_1, \cdots, \vv_{k}$, and let 
$q_0$ be the orthogonal projection of $q$ on $H_k$. Let $L$ be the affine subspace  of $\EE$ that passes through  $q_0$ and that is orthogonal to $H_k$.
Thus $H_k$ is  the affine subspace  of $\EE$ that passes through  $q_0$ and that is orthogonal to $L$. We may write these as
\begin{eqnarray}\label{perp}
(i) ~L = H_k^{\perp (q_0)}, ~~ (ii)~ H_k = L^{\perp (q_0)}.
\end{eqnarray}
It follows from
$$\dim L = \dim \EE - \dim H_k = (k+1)-(k-1) = 2$$ that $L$ is a 2-dimensional plane. 
Let $\Gamma$ be the circle in $L$ centered at $q_0$ and passing through $q$. If $t$ is any point on $\Gamma$, and if $1 \le j \le k$, then by Pythagoras' theorem, we have
$$\|t - \vv_j\|^2 = \|t - q_0 \|^2 + \|q_0 - \vv_j\|^2 = \|q - q_0 \|^2 + \|q_0 - \vv_j\|^2 = \|q - \vv_j\|^2,$$
and therefore 
\begin{eqnarray}\label{1jk}
\| t - \vv_j \| &=& \|q- \vv_j\|, \mbox{~(i.e., $\phi_j (t) = \phi_j (q)$) for $1 \le j \le k$}.
\end{eqnarray}
Also, the projection of $\vv_{k+1}$ on $L$ is not the center $q_0$ of $\Gamma$. Otherwise, $\vv_{k+1} q_0$ would be perpendicular to $L$ and hence $\vv_{k+1}$   would belong to $L^{\perp (q_0)}$, i.e., to $H_k$, contradicting the assumption that $\vv_1, \cdots, \vv_{k+1}$ are affinely independent. 
By Lemma \ref{III78},   the set $$W = \{\|\vv_{k+1} - t\| : t \in \Gamma\}$$ is an infinite set.

Let $r_j$, $1\le j\le k$, be defined by $r_j = \phi_j (q) = \| q - \vv_j \|$, and let $f(x) \in \RR[x]$ be defined by
$$f(x) = g(r_1,\cdots,r_k,x) \in \RR[x].$$
Then
$$f(x) = \sum_{j=0}^{d} g_j(r_1,\cdots,r_k) x^j.$$
This is a non-zero polynomial in $x$, since $g_d (r_1,\cdots,r_k) \ne 0$, by the choice of $q$; see (\ref{STAR}).
Also,  $f(w) = 0$ for all $w \in W$. In fact, if $w \in W$, then $w = \phi_{k+1} (t)$ for some $t$ in $\Gamma$,  and therefore
$$f(w) = g(\phi_1 (q), \cdots, \phi_k (q), \phi_{k+1} (t)) = 
g(\phi_1 (t), \cdots, \phi_k (t), \phi_{k+1} (t)),$$
by (\ref{1jk}). Hence $f(w)=0$, by (\ref{g-phi}). Thus $f(w) = 0$ for all $w$ in the infinite set $W$, a contradiction.
\hfill $\Box$

\bigskip
The next corollary is what we actually need in the proof of the main theorem, Theorem \ref{main} below. 

\begin{cor} \label{zzzz}
Let $n \ge 1$. Let $\EE$ be a Euclidean space of dimension  greater than or equal to  $n$, and let $\vv_1,\cdots,\vv_{n}$ be affinely independent points in $\EE$. Let $\EE_0$ be any everywhere dense subset of $\EE$, and let the functions $\phi_1,\cdots,\phi_n : \EE_0 \to \RR$ be defined by
$$\phi_j  (p) = \|p - \vv_j\|,~ 1 \le j \le n.$$
Then the functions $\phi_j$, $1 \le j \le n$, thought of as elements in the ring  of all real-valued functions on $\EE_0$,  are algebraically independent (over $\RR$). This is true in particular if $\EE \setminus \EE_0$ is a finite set. 
\end{cor}

\vspace{.1cm} \noindent {\it Proof.}
Let  $g \in \RR[T_1,\cdots,T_{n}]$ be a non-zero polynomial for which
$$g(\phi_1 (p), \cdots, \phi_n (p)) =0$$ for all $p\in \EE_0$. Since 
$g(\phi_1,\cdots,\phi_n)$ is continuous, its set of zeros is closed, and hence contains the closure $\EE$ of $\EE_0$. Thus
$$g\left(\|p-\vv_1\|,  \cdots, \|p - \vv_n\|\right) =0$$ for all $p\in \EE$.  \hfill $\Box$

\section{Answering  Question 1.1: The Main Result} \label{MAIN}
In this section, we establish, in Theorem \ref{main}  an  affirmative answer to Question \ref{Q1} posed in Section \ref{int}. Thus we prove  that the ideal $I$ of $R$ is indeed the principal ideal generated by $F$.


\begin{thm} \label{main}
Let $S = [\vv_1, \cdots, \vv_{d+1}]$, $d \ge 2$, be a regular $d$-simplex in $\RR^{d}$ having edge length $a$. Let $t_1, \cdots, t_{d+1}$ be the distances from an arbitrary point $\ttt$ in $\RR^d$ to the vertices of $S$. Let $R = \RR [T_1,\cdots,T_{d+1}]$  be the ring of polynomials over $\RR$ in the indeterminates $T_1, \cdots, T_{d+1}$, and let $I$ be the ideal of $R$ defined by
\begin{eqnarray}
I&=& \{ f \in R : f\left(t_1,\cdots, t_{d+1}\right) = 0 ~~\forall~ \ttt \in \RR^{d}\}. \label{II}
\end{eqnarray}
Let  \begin{eqnarray}
F &=& (d+1) \left( a^4 + \sum_{j=1}^{d+1} T_j^4\right) - \left( a^2 + \sum_{j=1}^{d+1} T_j^2\right)^2. \label{FF}
\end{eqnarray}
Then  $I$  is the principal ideal generated by $F$.  
\end{thm}

\vspace{.2cm} \noindent {\it Proof.} Let  $U = \RR^d \setminus \{ \vv_1,\cdots,\vv_{d+1}\}$.
Let $\AAA$ be the set of  (real-valued) functions  that are real analytic on $U$.
By Corollary \ref{xxx}, $\AAA$ is an integral domain. By Corollary \ref{yyy}, the functions $\phi_j : U \to \RR$, $1 \le j \le d+1$, defined by $\phi_j (\ttt) = \|\ttt - \vv_j\|$ belong to $\AAA$. 

Let $\Phi : \RR[T_1,\cdots,T_{d+1}]  \to \AAA$ be the ring $\RR$-homomorphism defined by $\Phi (T_j) = \phi_j$  for $1 \le j \le d+1$. The kernel of $\Phi$ consists of all $f \in R$ such that $f(t_1,\cdots,t_{d+1}) = 0$ for all $\ttt$ in $U$. By continuity, this is the set of all $f \in R$ such that $f(t_1,\cdots,t_{d+1}) = 0$ for all $\ttt$ in $\RR^d$, i.e., the ideal $I$.
The  image $\AAA_0$ of $\Phi$,  being a subring of the integral domain $\AAA$, is itself an integral domain. Since
$R/I \cong \AAA_0,$ by the first isomorphism theorem, it follows that $I$ is a prime ideal of $R$.

Also,  $\AAA_0$ contains the $d$  elements $\phi_1, \cdots, \phi_d$, and these are algebraically independent over $\RR$, by Corollary \ref{zzzz}. It follows from the definition that  $\td_{\RR} (QF(\AAA_0)) \ge d$, where $QF(.)$ denotes the quotient field.  Since $\AAA_0$ and $R/I$ are isomorphic as $\RR$-algebras, it follows that  $\td_{\RR} (QF(R/I)) \ge d$. Also, the integral domain $R/I$ is an affine $\RR$-algebra. Therefore $\td_{\RR} (QF(R/I)) = \dim (R/I)$. Therefore 
\begin{eqnarray*}
d &\le& \td_{\RR} (QF(\AAA_0)) ~=~\td_{\RR} (QF(R/I)) ~=~ \dim (R/I)\\
&\le& \dim (R) - \htt (I) ~=~ d+1 - \htt (I).
\end{eqnarray*}
Hence $\htt (I) \le 1$.
Since $I$ contains the polynomial $F$, which is irreducible by taking $n=t=d+1$ in Theorem \ref{Van}, it follows that $I$ contains  the prime ideal generated by $F$. Hence $\htt (I) = 1$. By Theorem \ref{ht=1}, $I$ is generated by $F$, as claimed.  \hfill $\Box$

\begin{rem}
{\em
When talking about a $d$-simplex, one usually assumes that $d \ge 2$, since a 1-simplex is a line segment with a poor geometry. However, it is legitimate to wonder whether Theorem \ref{main} still holds when $d=1$, and it may be interseting to know that it does not. In fact, if one takes $d=1$ and if one defines $I$ and $F$ as in Theorem \ref{main}, then it turns out that $F$ is not irreducible any more, as it factors into 
\begin{eqnarray*}
F&=& 2(a^4+T_1^4+T_2^4)-(a^2+T_1^2+T_2^2)^2\\
&=&(T_1+T_2+a)(T_1+T_2-a)(T_1-T_2+a)(T_1-T_2-a),
\end{eqnarray*} 
and that $I$ is the principal ideal generated by
\begin{eqnarray*}
F_0&=& (T_1+T_2-a)(T_1-T_2+a)(T_1-T_2-a),
\end{eqnarray*} 
and not by $F$. To see this, let $\vv_1, \vv_2$ be two distinct points in $\RR$, and let
\begin{eqnarray*}
I&=& \{f(T_1,T_2) \in \RR[T_1,T_2] : f(t_1,t_2) = 0 ~~\forall~ \ttt \in \RR\}.
\end{eqnarray*} 
Assuming that $\vv_1 < \vv_2$, we let $I_1$, $I_2$, and $I_3$ be the ideals defined by
\begin{eqnarray*}
I_1&=& \{f(T_1,T_2) \in \RR[T_1,T_2] : f(t_1,t_2) = 0 ~~\forall~ \ttt \le \vv_2\},\\
I_2&=& \{f(T_1,T_2) \in \RR[T_1,T_2] : f(t_1,t_2) = 0 ~~\forall~ \ttt \ge \vv_1\},\\
I_3&=& \{f(T_1,T_2) \in \RR[T_1,T_2] : f(t_1,t_2) = 0 ~~\forall~ \ttt \in (\vv_1, \vv_2)\}.
\end{eqnarray*} 
It is easy to see that the polynomials $$H_1=T_1-T_2-a \in I_1,~
H_2=T_1-T_2+a \in I_2,~H_3=T_1+T_2-a \in I_3.$$
Being of total degree 1, $H_i$, $1 \le i \le 3$,
is irreducible, and hence generates $I_i$. 
Thus if $H \in I$, then $H \in I_1 \cap I_2 \cap I_3$, and therefore $H_i$ divides $H$ for $1 \le i \le 3$. Since $H_1$, $H_2$, and $H_3$ are pairwise relatively prime, and since $\RR[T_1,T_2]$ is a unique factorization domain, it follows that $H_1H_2H_3$ divides $H$. This shows that $I$ is generated by $H_1H_2H_3$, as claimed.
}
\end{rem}

\section{Questions for further research} \label{Q}

\begin{que} \label{Q2}
{\em Let $S = [\vv_1,\cdots,\vv_{d+1}]$ be a regular $d$-simplex of side length $t_0$ in $\RR^d$. For every $\ttt \in \RR^d$,  let $t_j$ be the distance from $\ttt$ to $\vv_j$, $1\le j \le d+1$. Let 
$R_0 = \RR[T_0,T_1,\cdots,T_{d+1}]$ 
be the ring of polynomials in the indeterminates   $T_0$, $T_1$, $\cdots$,  $T_{d+1}$ over the field $\RR$ of real numbers. Allowing $t_0$ to vary, 
let 
$I_0$  be the ideal in $R_0$ defined by
\begin{eqnarray}
I_0&=& \{ f \in R_0 : f\left(t_0,t_1,\cdots, t_{d+1}\right) = 0 ~~\forall~ \ttt \in \RR^{d}\}, \label{I0}
\end{eqnarray}
and let  \begin{eqnarray}
F_0 &=& (d+1) \left( T_0^4 + \sum_{j=1}^{d+1} T_j^4\right) - \left( T_0^2 + \sum_{j=1}^{d+1} T_j^2\right)^2. \label{F0}
\end{eqnarray}
Is $I_0$  generated by $F_0$?}
\end{que}


\begin{que} \label{Q3}
{\em Suppose that $S=[\vv_1,\cdots,\vv_{d+1}]$ is a regular $d$-simplex of side length $a$ in $\RR^d$, and suppose that the positive numbers
$t_1, \cdots, t_{d+1}$ satisfy (\ref{rel}). Does there exist a point 
 $\ttt \in \RR^d$ such that  the distance from $\ttt$ to $\vv_j$, $1\le j \le d+1$, is $t_j$?}
 \end{que}	

\begin{que} \label{Q4}
{\em Let $S = [\vv_1,\cdots,\vv_{d+1}]$ be a regular $d$-simplex of side length $t_0$ in $\RR^d$, and let $\Gamma$ be the circumsphere of $S$.  For every $\ttt \in \Gamma$,  let $t_j$ be the distance from $\ttt$ to $\vv_j$, $1\le j \le d+1$. Let $R = \RR[T_1,\cdots,T_{d+1}]$ 
be the ring of polynomials in the indeterminates   $T_1$, $\cdots$,  $T_{d+1}$ over the field $\RR$ of real numbers, and let 
$J$  be the ideal in $R$ defined by
\begin{eqnarray}
J&=& \{ f \in R : f\left(t_1,\cdots, t_{d+1}\right) = 0 ~~\forall~ \ttt \in \Gamma\}. \label{J}
\end{eqnarray}
Clearly $J$ contains the polynomial $F$ defined by
\begin{eqnarray*}
F &=& (d+1) \left( a^4 + \sum_{j=1}^{d+1} T_j^4\right) - \left( a^2 + \sum_{j=1}^{d+1} T_j^2\right)^2.
\end{eqnarray*}
It is also known, and not difficult to prove, that $J$ contains the polynomials $G$ and $H$ defined by
\begin{eqnarray*}
G&=& \left(  \sum_{j=1}^{d+1} T_j^2\right) - d a^2, ~~~H ~=~ \left(  \sum_{j=1}^{d+1} T_j^4\right) - d a^4.
\end{eqnarray*}
Thus one may ask whether $J$ is generated by $F$, $G$, and $H$. However, it is easy to check that
\begin{eqnarray}
(d+1)H&=&F+((d+1)a^2+G)^2 - (d+1)^2 a^4,
\end{eqnarray}
and hence $F$ is generated by $G$ and $H$, and $H$ is generated by $F$ and $G$. Thus we ask
\begin{quote}\begin{quote}
	 Is $J$ generated by $G$ and $H$?  Is $J$ generated by $F$ and $G$?
\end{quote}\end{quote}
}
\end{que}

\begin{que} \label{Q5}
{\em Is the Soddy relation essentially the only one? Given positive numbers that satisfy the Soddy relation, do there exist spheres having these numbers as radii and  mutually touching each other?}
\end{que}

\begin{que} \label{Q6}
{\em If, instead of taking a regular $n$-simplex, one takes a general $n$-simplex with given edge lengths, then the relation among the distances of an arbitrary point in its affine hull is expected to be complicated. In face, this problem is addressed for a triangle in \cite{Bottema}, and the relation is found to be quite unmanageable. However, one may try to consider a tetrahedron which is not quite general. For example, a reasonable analog of the triangle in 3-space is, besides the regular tetrahedron, the tetrahedron having congruent faces. These tetrahedra are called {\it equifacial}, and has attracted much attention.}
\end{que}

\vspace{.3cm} \noindent {\bf Acknowledgment.} The first-named author is supported by Research Grant 2014/15 from Yarmouk University.

\end{document}